\documentclass{article}

\usepackage[american]{babel}
\usepackage{amsfonts,amsmath,amssymb,epsf,epsfig}
\usepackage{xypic}
\xyoption{all}
\input{xypic}
\newdir{ >}{{}*!/-8pt/\dir{>}}

\def\R{\mathbb{R}}

\def\C{\mathbb{C}}

\def\Z{\mathbb{Z}}

\def\P{\mathbb{P}}

\def\Spec{{\rm Spec}}
\def\Hilb{{\rm Hilb}}
\def\Vect{{\rm Vect}}
\def\Der{{\rm Der}}
\def\Hom{{\rm Hom}}
\def\id{\rm id}

\def\pr{\noindent $\bf{Proof.}$\quad}     
\def\fin{\hfill$\square$\\}           
\newtheorem{theo}{Theorem}
\newtheorem{defi}{Definition}
\newtheorem{rem}{Remark}
\newtheorem{prop}{Proposition}
\newtheorem{cor}{Corollary}

\newtheorem{lem}{Lemma}
\begin{document}
\title{Deformations of Lie algebras of vector fields arising from families of schemes} 

\author{Friedrich Wagemann \\
        Universit\'e de Nantes}
   
\maketitle


\section*{Introduction}

The goal of the present paper is to construct examples of global deformations 
of vector field Lie algebras in a conceptual way. Fialowski and Schlichermaier 
\cite{FiaSch}
constructed global deformations of the infinitesimally and formally rigid Lie 
algebra of polynomial vector fields on the circle and of its central extension,
the Virasoro algebra. The Lie algebra of polynomial vector fields on the circle
is here replaced by/seen as the Lie algebra of meromorphic or rational vector 
fields on the Riemann sphere admitting poles only in the points $0$ and 
$\infty$. In this context one gets non trivial 
deformations from an affine family of curves by first deforming it as 
a projective family with marked points and then extracting the points.
Fialowski and Schlichermaier get in this way non-trivial deformations of Lie 
algebras, and the underlying families of curves present singularities. 

In the attempt of producing deformations of vector field Lie algebras from 
deformations of the underlying pointed algebraic variety in a general 
framework, we are led to a notion of {\it global deformation} which is 
different from the one used by Fialowski and Schlichenmaier, one 
which is closer to deformations in algebraic geometry. A first goal 
is to compare these global deformations with the corresponding notion 
from Fialowski-Schlichermaier's article.

We then stick to the notion of deformations of the Lie algebra of vector fields
on a pointed algebraic curve imposed by the deformation of the underlying 
curve, and show its close relation to the moduli space of pointed curves. 
In order to formulate this relation, we show that the ``space of deformations''
carries the geometric structure of a $\C$-stack. This is a kind of functor
from the category of model spaces, here affine schemes, to the category of
groupoids (in order to capture the idea of considering deformations up to 
isomorphism). It has to satisfy conditions to provide construction
of geometric objects by gluing local data, making it ressemble the functor
of points of a scheme. The link between the moduli stack ${\cal M}_{g,n}$ 
and the deformation stack ${\cal D}{\rm ef}$ is a morphism of stacks $I$.
It realizes a family of marked projective curves as a deformation of the Lie 
algebra of regular vector fields on the affine curve obtained from extracting 
the points, by extracting the divisor associated to the union of the marked 
points, and then taking the Lie algebra of sections of the relative tangent 
bundles of the family. 

We show furthermore that $I$ is almost a monomorphism. The idea that this 
might be true comes from Pursell-Shanks' theory of describing the underlying 
manifold (and more geometric objects related to it) by its Lie algebra of 
tangent vector fields. We generalize this theory to a relative (affine) 
setting, i.e. we show that over $\C$, $\Der_A(B)\cong\Der_A(B')$ implies an 
$A$-isomorphism $B\cong B'$.  

Let us observe that all which is discussed in this paper can be easily 
transposed to current algebras, instead of vector field Lie algebras. 
We hope that 
the stack approach to deformations is useful from the two points of view:
our believe is that on the one hand it may give rise to a cohomology theory
gouverning global deformations by associating some standard homological algebra
tool to the above morphism of stacks $I$ (its cotangent complex ?). 
For this direction, it might be necessary
to show first that some substack of our deformation stack is an algebraic 
stack, cf the remark at the end of section $5$.
 
\noindent{\bf Acknowledgements:} the author thanks Alice Fialowski for the 
discussion at Strasbourgh which initiated the present study. He thanks 
furthermore Pierre-Emmanuel Chaput for help with lemma $3$ and other useful 
remarks, as well as Matthias Borer for useful discussions. Some of this work
was discussed with Martin Schlichenmaier at Luxembourg, and my thanks go
to him as well.

\section{Global deformations of general and of product type}

Let $\C$ be the base field; we are conscious that some results are valid in a 
much more general framework.

In this first section, we compare two types of global deformations of Lie 
algebras. On the one hand, we have global deformations as 
Fialowski-Schlichenmaier \cite{FiaSch}
define them. We will call these {\it of product type}. In the following, 
$A$ is an associative commutative unital $\C$-algebra (and all $\C$-algebras
will be supposed associative commutative and unital).

\begin{defi}
For a Lie algebra ${\mathfrak g}_0$, a global deformation (of product type) 
with base $A$ is an $A$-Lie algebra structure $[-,-]_{\lambda}$ on the tensor 
product ${\mathfrak g}_0\otimes_{\C} A$, together with a morphism of Lie 
algebras $\id\otimes\epsilon:{\mathfrak g}_0\otimes A\to 
{\mathfrak g}_0\otimes \C\simeq{\mathfrak g}_0$, where $\epsilon$ is the 
augmentation of the (augmented) algebra $A$.
\end{defi}

Concretely, the bracket $[-,-]_{\lambda}$ is antisymmetric, satisfies the 
Jacobi identity, and the following two identities for $a,b\in A$ and $x,y\in
{\mathfrak g}_0$ owing to $A$-linearity of the bracket and the 
fact that $\epsilon\otimes\id$ is supposed to be a Lie algebra homomorphism:

\begin{itemize}
\item[(1)] $[x\otimes a,y\otimes b]_{\lambda}\,=\,[x\otimes 1,y\otimes 
1]_{\lambda}\,\,(ab)$
\item[(2)] $({\id}\otimes\epsilon)([x\otimes 1,y\otimes 1]_{\lambda})\,=\,
[x,y]\otimes 1$.
\end{itemize}

By condition $(1)$, it is enough to give the commutators 
$[x\otimes 1,y\otimes 1]_{\lambda}$ for all $x,y\in{\mathfrak g}_0$ 
in order to describe the deformation. By condition $(2)$, their bracket 
has the form
$$[x\otimes 1,y\otimes 1]_{\lambda}\,=\,[x\otimes 1,y\otimes 1]\,+
\,\sum_iz_i\otimes a_i.$$

We call the notion with which we want to contrast these global deformations of
product type {\it general global deformations}. We take over the term 
``global'' from earlier work on the subject, notwithstanding the fact that
these deformations are ``only'' on affine open sets (namely, on the affine open
set given by $\Spec(A)$, the spectrum of $A$).

\begin{defi}
Let ${\mathfrak g}_0$ be a $\C$-Lie algebra. A general global deformation 
of ${\mathfrak g}_0$ is an $A$-Lie algebra ${\mathfrak g}$ together with 
an augmentation $\epsilon:A\to \C$ such that
${\mathfrak g}\otimes_A \C\simeq{\mathfrak g}_0$ as $\C$-Lie algebras
(where $A$ acts on $\C$ via $\epsilon$).
\end{defi}

One detail is hidden in this definition. ${\mathfrak g}\otimes_A \C$ is a 
$\C$-Lie algebra: indeed, ${\mathfrak g}\otimes_A \C$ is a 
$\C$-Lie algebra with the current algebra bracket 
$[x\otimes \alpha,y\otimes \beta]:=[x,y]\otimes \alpha\beta$ which we call 
the {\it trivial bracket} here. Then the space
$$I:={\rm span}_{\C}\{ax\otimes\alpha\,-\,x\otimes\epsilon(a)\alpha\,|\,
x\in{\mathfrak g}_0,\,\,a\in A,\,\,\alpha\in \C\,\}$$
is a Lie-ideal as follows directly from the $A$-linearity of the bracket 
on ${\mathfrak g}$, so the quotient ${\mathfrak g}\otimes_A \C$ is a 
Lie algebra.

\begin{prop}
A deformation of product type $({\mathfrak g}_0\otimes_{\C} A,[,]_{\lambda},
\epsilon)$ is a general global deformation.
\end{prop}

\pr It is clear that $({\mathfrak g}_0\otimes_{\C} A)\otimes_A \C\simeq 
{\mathfrak g}_0$ as $\C$-modules. The morphisms can be given by
$$f:({\mathfrak g}_0\otimes_{\C} A)\otimes_A {\C}\to {\mathfrak g}_0,\,\,
\,\,\,f(x\otimes a\otimes\alpha)=\alpha\epsilon(a)x\,=\,
\alpha({\id}\otimes\epsilon)(x\otimes a)$$
and 
$$g:{\mathfrak g}_0\to({\mathfrak g}_0\otimes_{\C} A)\otimes_A \C,\,\,
\,\,\,g(x)=x\otimes 1\otimes 1.$$
Using that the $A$-module structure on $\C$ is given by the augmentation 
$\epsilon$, it is easy to see that $f$ and $g$ are mutually inverse. 
Let us show that $f$ is a morphism of Lie algebras:

\begin{eqnarray*}
f([x\otimes a\otimes\alpha,y\otimes b\otimes\beta])&=&f([x\otimes a,y\otimes 
b]_{\lambda}\otimes\alpha\beta) \\
&=&f([x\otimes 1,y\otimes 1]_{\lambda}\,(ab)\,\otimes\alpha\beta) \\
&=&f([x\otimes 1,y\otimes 1]_{\lambda}\otimes\epsilon(ab)\alpha\beta) \\
&=&\epsilon(ab)(\alpha\beta)\,\,({\id}\otimes\epsilon)([x\otimes 1,y\otimes 
1]_{\lambda}) \\ &=&\epsilon(ab)(\alpha\beta)\,\,[x,y]. 
\end{eqnarray*}

where we used conditions $(1)$ and $(2)$, and on the other hand

$$[f(x\otimes a\otimes\alpha),f(y\otimes b\otimes\beta)]\,=\,\epsilon(ab)
\alpha\beta\,\,[x,y].$$
This ends the proof of the lemma.\fin

Thus, the notion of general global deformations includes the deformations of 
product type, as it should. The reason for considering these more general 
deformations is that in general the Lie algebra deformations of a Lie algebra 
of vector fields coming from deformations of the underlying variety are not of 
product type.  

\section{A map from families of schemes to deformations of Lie algebras}

A(n affine) {\it deformation} of a scheme $X_0$ over $\C$ is a scheme 
$X$ over a $\C$-algebra $A$ and a closed point $0\in \Spec(A)$ (with residue
field $\C$) such that the fiber $(X)_0$ of $X$ in $0$ satisfies $(X)_0\simeq 
X_0$ as $\C$-schemes. Let us denote by $\pi:X\to \Spec(A)$ the structure 
morphism. Usually we suppose furthermore that $\pi:X\to \Spec(A)$ is flat,
quasi-compact, surjective and finitely presented. Recall that when 
$\pi$ is of finite type and $A$ noetherian, $\pi$ is automatically of 
finite presentation. As one encounters non-noetherian bases in descent theory,
finite presentation is the good hypothesis to impose (for example, in the 
realm of fppf topology). We then 
regard $X$ as a {\it family of schemes} deforming $X_0$. By definition of 
the fiber, $(X)_0\simeq X_0$ means $X\times_{\Spec(A)}\Spec(\C)\simeq X_0$.

In order to get started, let us suppose that the family is {\it affine}, 
i.e. $X=\Spec(B)$, and that $X_0$ is affine also, i.e. $X_0=\Spec(B_0)$. 
We are interested in the Lie algebra of vector fields on $X$ which are 
tangent to the fibers of $\pi$, i.e. in the Lie algebra of $A$-linear 
derivations of $B$, denoted $\Der_A(B)$.

The following lemma expresses the behavior of $A$-linear 
derivations of $B$ under base change.

\begin{lem}
Let $A$ be a $\C$-algebra, and $A'$ and $B$ be $A$-algebras. 
Let $B':=B\otimes_A A'$ the coproduct in this situation. 
Suppose that $B'$ is a flat $B$-module and that the $B$-module 
of K\"ahler differentials $\Omega_A(B)$ has a finite presentation. Then
$$\Der_{A'}(B')\,\simeq\,\Der_A(B)\otimes_A A'$$
as $A'$-modules.
\end{lem}

\begin{rem}
Translating into scheme language, i.e. $X=\Spec(A)$, $X'=\Spec(A')$, 
$Y=\Spec(B)$, $Y'=\Spec(B')$ is the fibered product $Y'=X'\times_X Y$, 
and the hypotheses are satisfied in case the induced morphisms 
$\pi:Y\to X$ and $X'\to X$ are of finite presentation and flat. 
This situation is called base change (in the fppf topology).
\end{rem} 

\pr This follows in a straight forward manner from the corresponding base 
change for the modules of K\"ahler differentials. Here is an explicit proof: 

The $A'$-module of K\"ahler differentials $\Omega_{B'/A'}$ of the 
$A'$-algebra $B'$ represents the functor $\Der_{A'}(B',-)$ which 
implies in particular an isomorphism of $B'$-modules
$$\Der_{A'}(B')\simeq \Hom_{B'}(\Omega_{B'/A'},B').$$
The hypothesis $B\otimes_A A'= B'$ implies 
$\Omega_{B'/A'}\simeq\Omega_{B/A}\otimes_BB'$, which in turn gives
$$\Hom_{B'}(\Omega_{B'/A'},B')\simeq\Hom_{B'}(\Omega_{B/A}\otimes_B B',B').$$
Then we apply the ``change of rings'' isomorphism (meaning the adjointness of
$-\otimes_BB'$ and $\Hom_{B'}(B',-)$) 
$$\Hom_{B'}(\Omega_{B/A}\otimes_B B',B')\simeq \Hom_B(\Omega_{B/A},B').$$
In order to retranslate the result into $\Der_A(B)\otimes_BB'$, we need 
for a $B$-module $M$ a natural isomorphism
$$\Hom_B(M,B')\simeq\Hom_B(M,B)\otimes_BB'.$$
There is obviously such an isomorphism for $M=B$. Furthermore, the 
claim holds true for a $B$-module $M$ of finite presentation as is 
easily checked ($B'$ is a flat $B$-module by hypothesis). But 
$\Omega_{B/A}$ has a finite $B$-presentation (here we also use the 
noetherian hypothesis). Thus we can conclude 
$$\Hom_B(\Omega_{B/A},B')\simeq\Hom_B(\Omega_{B/A},B)\otimes_BB',$$
and
$$\Hom_B(\Omega_{B/A},B)\otimes_BB'\simeq \Der_A(B)\otimes_BB'.$$
Finally
$$\Der_A(B)\otimes_BB'=\Der_A(B)\otimes_B(B\otimes_A A')\simeq
\Der_A(B)\otimes_A A'.$$
\fin

\begin{rem}
It is easy to check that the above isomorphism
$$\Der_{A'}(B')\,\simeq\,\Der_A(B)\otimes_A A'$$
is an isomorphism of $A'$-Lie algebras when we give the right hand side 
the bracket $[x\otimes a,y\otimes b]:=[x,y]\otimes ab$. This follows from 
the fact that the bracket of derivations $D$, $D'$ is on both sides 
given by $[D,D']=D\circ D' - D'\circ D$.
\end{rem} 

\begin{lem}
Let $A$ be a $\C$-algebra, $B\to A$ be an $A$-algebra, and $A\to\C$ be a 
$\C$-point of the affine scheme $\Spec(A)$. Suppose that the $B$-module 
$\Omega_{B/A}$ is projective. Then there is an isomorphism of 
$\C$-modules
$$\Der_A(B)\otimes_A \C\simeq \Der_{\C}(B_0).$$
\end{lem}

\pr This is a slightly different situation from that in lemma $3$ where 
the module of derivations of the fiber (product) is just the tensor product
of the initial module of derivations.

Its justification stems from the fact that in the above proof of lemma $3$, 
the natural isomorphism
$$\Hom_B(M,B')\simeq\Hom_B(M,B)\otimes_BB'$$
holds not only for $M=B^n$, but also for a direct factor of $B^n$, i.e. for
each projective $B$-module, by additivity of the concerned functors. 
Observe that we do not need any flatness assumption here.
\fin 

\begin{rem}
Translated into scheme language, this lemma shows that the Lie algebra of 
derivations of a fiber of $\pi:\Spec(B)\to\Spec(A)$ is the tensor product 
of $\Der_A(B)$ by $\C$, seen as an $A$-module via the augmentation $A\to\C$, 
in case $A$ is smooth. More generally, this is true at each smooth point of 
$A$. 

Translated into the language of deformations exposed in section $1$, it shows 
that deformations of the underlying affine scheme $X_0=\Spec(B_0)$ give rise 
to general global deformations of the Lie algebra of derivations $\Der(B_0)$.
\end{rem}

In conclusion, there is a map, which we call $I$, assigning to an affine 
family of schemes $X$ (over $A$) deforming a given affine scheme $X_0$ a 
general global $A$-deformation of the Lie algebra $Vect(X_0)$. 

\begin{rem}
Note that the map $I$ does not define a functor in general, 
as $\Der_A(B)$ is not functorial in $B$; the point is that pushforward or 
pullback of vector fields along morphisms works for isomorphisms.

For example, one can define a genuine functor when restricting the source 
category to a groupoid of some kind of varieties with allowed morphisms 
being only the isomorphisms. 
\end{rem}

\begin{cor}
In the case of a product $Y'=Y\times X'$ and $X=\Spec(\C)$, i.e. a 
trivial deformation, one gets an isomorphism of Lie algebras 
$$\Der_{A'}(B')\,\simeq\,\Der_{\C}(B)\otimes_{\C} A'.$$
\end{cor}

\begin{rem}
Here the right hand side carries the current algebra bracket 
$[x\otimes a,y\otimes b]:=[x,y]\otimes ab$ which is considered as the 
trivial deformation, thus we see that trivial deformations of the 
underlying variety give rise to trivial global deformations of product type.
\end{rem} 

\begin{rem}
The difference between the notions of general global deformations and 
of global deformations of product type is that the latter considers only 
global deformations which are, as $A$-modules, tensor products with $A$ over
$\C$. But these are then free $A$-modules, whereas in general the Lie 
algebra of $A$-linear derivations on an algebra $B$ has no reason to be a 
free $A$-module. It would be interesting to have a clear criterion when
this is the case.

Nevertheless, in case the special point $0\in \Spec(A)$ is a generic point 
in the moduli space ${\mathcal M}_{g,n}$, $\Der_A(B)$ is a free $A$-module,
because in this situation, the sheaf of Krichever-Novikov algebras is locally 
free, and the local basis trivializes the module, see \cite{Sch} p. 743. 
We owe this remark to M. Schlichenmaier.
\end{rem}

\section{Extracting affine families from proper families}

Let $S$ be some scheme over $\C$ and $n,g$ integers such that $2g-2+n>0$. 
A proper flat (surjective) family $\pi:C\to S$  
whose geometric fibers $C_s$ are reduced
connected curves of genus $g={\rm dim}\,\,H^1(C_s,{\cal O}_{C_s})$ 
with at most ordinary double points is
called a {\it stable curve} of genus $g$ over $S$. An {\it $n$-pointed stable 
curve} is a stable curve of genus $g$ over $S$ together with $n$ sections 
$\sigma_1,\ldots,\sigma_n:S\to C$ of $\pi$ such that the images 
$P_i=\sigma_i(s)$ are disjoint and smooth for all $i=1,\ldots,n$ and all 
$s\in S$.
Furthermore, the number of points where a nonsingular rational 
component $E$ of $C_s$ meets the rest of $C_s$ plus the number of points 
$P_i$ on $E$ must be at least $3$, cf \cite{Knut} definition $1.1$ p.162.
In the same vein, we call an {\it affine stable curve} a flat (surjective)
family $\pi:C\to S$ whose geometric fibers are affine stable reduced and 
connected curves of genus $g$ with at most ordinary double points. All our
affine stable curves will be obtained from $n$-pointed stable curves by 
extracting the marked points and one may take this for the definition of an
affine stable curve. We will then be in position to apply the 
results of section $2$: 

\begin{lem}
Let $\pi:C\to S$ be a punctured stable curve, and $\sigma_1,\ldots,\sigma_n$ 
the corresponding sections of $\pi$.  
Denote by $\triangle$ the union $\triangle=\bigcup_{i=1}^n\sigma_i(S)$. 
Then $\pi:C\setminus\triangle\to S$ is an affine stable curve and an 
affine morphism.
\end{lem}

\pr It is shown in \cite{Knut} p. $173$ corollary $1.9$ that some tensor
power of the sheaf $\omega_{C_s/\C}(P_1+\ldots+P_n)$ is very ample, 
where $\omega_{C_s/\C}$ is the dualizing sheaf of $C_s$ over $\C$.
Here we denote by $C_s=\pi^{-1}(s)$ for $s\in S$, and $P_i=\sigma_i(s)$ for
$i=1,\ldots,n$. Let us denote by ${\cal L}_{\triangle}$
the line bundle $\omega_{C/S}(\triangle)$ on $C$. 
By Grauert's theorem (cf \cite{Hart} pp. 288--289), we have 
$R^i\pi_*({\cal L}^{\otimes n}_{\triangle})=0$ for all $i>0$ and $n>>0$. 
The local to global spectral sequence gives therefore 
$H^i(C,{\cal L}^{\otimes n}_{\triangle})=0$ for all $i>0$ and $n>>0$. 
${\cal L}_{\triangle}$ is then an ample line bundle on $C$, and we use 
some power ${\mathcal L}^{\otimes k}_{\triangle}$ of ${\mathcal L}_{\triangle}$
to embed $C$ into projective space $\P^N_A$; let $i:C\to\P^N_A$ be 
corresponding embedding. Then $C\setminus\triangle$ maps under $i$ into the 
complement of a hyperplane, and $i|_{C\setminus\triangle}$ is affine. 
The claim follows.\fin 

Now we pass to moduli space; there are (at least) two points of view on the
moduli space of curves, and it is the necessity of the situation which
determines which one to take. One point of view is the following (the GIT
or coarse moduli space point of view): the {\it moduli space} 
(of stable connected projective algebraic curves over $\C$ of fixed genus 
$g$ with $n$ marked points) is a scheme whose geometric points represent 
the isomorphism classes of stable punctured curves. 
We observe in lemma $4$ below that isomorphic punctured curves 
give rise to isomorphic deformations, and get thus that $I$ factors to a 
map from the set of $S$-points of the moduli space to global deformations 
of vector field Lie algebras over $S$. 

More precisely, let $\pi:C\to S$ and $\pi':C'\to S$ be two isomorphic 
punctured curves by an 
$S$-morphism $\phi$ which respects the punctures. Let $U=\Spec(B)\subset C$ be 
an affine open set such that $\pi|_U$ and $\pi|_{\phi(U)}$ are affine families 
of schemes. Denote by $B$ and $B'$ the corresponding algebras, i.e. 
$U=\Spec(B)$ and $\phi(U)=\Spec(B')$. One may assume that 
$\pi(U)=\pi'(\phi(U))=\Spec(A)$. 

\begin{lem}
In the above notation, the isomorphic curves $\pi:C\to S$ and $\pi':C'\to S$ 
give rise to isomorphic Lie algebras $\Der_A(B)$ and $\Der_A(B')$.  
\end{lem}

On the other hand, let $\pi:C\to S$ and $\pi':C'\to S$ be two punctured curves,
let $T=\Spec(A)\subset S$ be an affine open subset of the parameter space and 
denote by $U=\Spec(B)$ and $U'=\Spec(B')$ the affine open sets $U=\pi^{-1}(T)$ 
and $U'=(\pi')^{-1}(T)$ corresponding to the $\triangle$ (resp. $\triangle'$)
complement in $C$ (resp. $C'$).

\begin{theo}
In the above notation, suppose 
\begin{itemize}
\item either that $B$ and $B'$ are integral, normal, of finite type over $A$ 
noetherian,
\item or that all fibers of $C$ and $C'$ are smooth curves,
\item or that $\Der_A(B)$ (resp. $\Der_A(B')$) is a free $B$ (resp. $B'$) 
module of rank $1$.
\end{itemize}
Then $\Der_A(B)\simeq\Der_A(B')$ as Lie algebras 
over $A$ implies that there exists a $T$-isomorphism $\phi:U\to U'$ which is 
compatible with the punctures.
\end{theo}

\begin{rem}
The special case where $T$ is a closed point of residue field $\C$ is easily 
solved: the claim is then the fact that an isomorphism between punctured 
projective curves $f:(X\setminus\{p_1,\ldots,p_n\})\simeq(Y\setminus\{q_1,
\ldots,q_n\})$ can be extended to an isomorphism $X\simeq Y$ respecting the 
punctures (use proposition 6.8 in \cite{Hart} p. 43 and induction applied to 
$i\circ f:(X\setminus\{p_1,\ldots,p_n\})\to(Y\setminus\{q_1,\ldots,q_n\})
\hookrightarrow Y$). 

The proof of the above theorem follows from Pursell-Shanks theory, and is 
given in the next section.
\end{rem}

One deduces from the theorem that we get a well-defined {\it injective} map 
from the set of $S$-points of the moduli space of punctured smooth curves 
(deforming a fixed punctured curve $(X_0,p_1,\ldots,p_n)$) to the set of 
isomorphism classes of deformations of the Lie algebra 
$Vect(X_0\setminus\{p_1,\ldots,p_n\})$ with base $S$.

\section{Pursell-Shanks theory}

In this section, we develop another key ingredient of our study of the 
moduli space of deformations of a given Lie algebra of regular vector fields 
on an affine curve. Here we recall how the Lie algebra of vector fields 
on a manifold encodes geometric objects of the manifold. This theory is 
due to Pursell-Shanks, Omori, Amemiya \cite{Amem}, Grabowski \cite{Grab}, 
and in our framework to Siebert \cite{Sieb}.

To give an example of the kind of proposition studied in this theory, 
Siebert proves that two complex normal reduced irreducible affine algebraic 
varieties $Y=\Spec(B)$ and $Y'=\Spec(B')$ over $\C$ are isomorphic if and 
only if their Lie algebras of vector fields are isomorphic as Lie algebras 
(over $\C$). The method of proof describes the points of the underlying 
variety in terms of data attached to the Lie algebra of vector fields 
${\cal L}$. Indeed, he proves for irreducible $Y$ a bijection between 
the regular points of $Y$ and the finite-codimensional maximal subalgebras 
$L$ of ${\cal L}$ such that $L$ does not contain a proper Lie ideal. 
Call the set of these subalgebras ${\mathcal M}$. The idea is that such 
an $L$ is the Lie subalgebra of vector fields vanishing at the point 
it represents. This bijection can be used to express the Zariski topology 
on $Y$ in Lie algebraic terms. Then Siebert defines an algebra of functions 
$B({\cal L})$ (cf definition $10$ \cite{Sieb}) on $Y$ by taking 
those functions $f:{\mathcal M}\to\C$ such that for all 
$\delta\in{\cal L}$, there is a $\theta\in{\cal L}$ such that 
$f(L)\delta - \theta\in L$ for all $L\in{\mathcal M}$. 
$B({\cal L})$ encodes thus the coefficient functions of 
the vector fields. Siebert shows that for ${\cal L}=\Der_{\C}(B)$ with 
$Y=\Spec(B)$ irreducible of finite type over $\C$, $B({\cal L})$ is an 
integral extension of $B$. Normality of $B$ then forces $B({\cal L})=B$. 
In this sense, Siebert does not construct an explicit isomorphism 
$Y\simeq Y'$, but he deduces the {\it existence} of an isomorphism from the
translation of the geometric structures into Lie algebraic objects. 

Our idea is that the above theory reflects that the map $I$ is injective in 
a certain sense. For this we need a {\it relative version} of 
the above statement, namely $\Der_A(B)\simeq\Der_A(B')$ if and only if 
$\Spec(B)\simeq\Spec(B')$ by an isomorphism over $\Spec(A)$. 
The following statement will be the main result of this section. 

\begin{prop}
Suppose 
\begin{itemize}
\item either that $B$ and $B'$ are integral, normal, of finite type over $A$ 
noetherian,
\item or that all fibers of $C$ and $C'$ are smooth (cf notations section $3$),
\item or that $\Der_A(B)$ (resp. $\Der_A(B')$) is a free $B$ (resp. $B'$) 
module of rank $1$.
\end{itemize}
Then $\Der_A(B)\simeq\Der_A(B')$ 
if and only if $\Spec(B)\simeq\Spec(B')$ by an isomorphism over $\Spec(A)$.
\end{prop}

\pr Denote as before $\Der_A(B)={\cal L}$. It is a modular Lie algebra in the
sense of definition $1$ of \cite{Sieb} (with respect to $B$). 
To every algebra homomorphism $\phi:A\to\C$, i.e. to every $\C$-point $x$ of
$\Spec(A)=X$, we associate the fiber $\pi^{-1}(x)$ of $\pi:Y=\Spec(B)\to X$. 
Now we apply Siebert's theory to describe the smooth fibers $\pi^{-1}(x)$:
$\phi$ gives rise to the $\C$-algebra $B\otimes_A\C$ where $A$ acts on $\C$
via $\phi$ and on $B$ via $A\to B$. For each point $\phi=x\in X$, there is a 
projection 
$$\pi_{\phi}:\Der_A(B)\to \Der_A(B)\otimes_A\C\to \Der_A(B\otimes_A\C)$$
given by the composition of $X\mapsto X\otimes 1$ and 
$$D\otimes 1\mapsto(b\otimes 1\mapsto D(b)\otimes 1).$$ 
In this way we reduce at a point $x$ the problem
to the fiber $\pi^{-1}(x)$, and Siebert's theory implies for $\pi^{-1}(x)$ 
smooth that
$$\pi^{-1}(x)\,=\,{\cal M}_{\phi},$$
the space of maximal, finite-codimension subalgebras $L_{\phi}$ of 
${\cal L}_{\phi}:=\Der_{\C}(B\otimes_A\C)$, and furthermore that the algebra 
of regular functions ${\rm Reg}(\pi^{-1}(x))$ is 
$${\rm Reg}(\pi^{-1}(x))\,=\,{\cal B}_{\phi}$$
with
$${\cal B}_{\phi}\,:=\,\{f\in\C^{{\cal M}_{\phi}}\,|\,\forall\delta\in
{\cal L}_{\phi}\,\exists\theta\in{\cal L}_{\phi}\,:\,f(L_{\phi})\delta -
\theta\in L_{\phi}\,\,\,\,\forall L_{\phi}\in{\cal M}_{\phi}\}.$$
In order to describe not only the fibers, but all points of $Y$, consider 
now the set of pairs
$${\cal M}\,=\,\{(\phi,L_{\phi})\,|\,\phi=x\in X,\,\,\,
{\rm and}\,\,L_{\phi}\in {\cal M}_{\phi}\}$$
and 
$${\cal B}\,=\,\{f\in\C^{\coprod_{\phi}{\cal M}_{\phi}}\,|\,\forall\delta\in
{\cal L}\,\exists\theta\in{\cal L}\,:\,\pi_{\phi}(f(L_{\phi})\delta -
\theta)\in L_{\phi}\,\,\,\,\forall L_{\phi}\in{\cal M}_{\phi}\,\,\,{\rm and}\,
\,\,\forall\phi=x\in X\}.$$
Observe that we first have to project the expression $f(L_{\phi})\delta -
\theta$ onto the fiber via $\pi_{\phi}$ in order to state that the result 
should vanish at $L_{\phi}$, seen as a point of ${\cal M}_{\phi}$.

Let us show that ${\cal B}$ is an integral extension of $B$. Indeed, clearly
$B\subset{\cal B}$, because ${\cal L}$ is a $B$-module. 

If $B$ is a finitely generated $A$-algebra, there exists an epimorphism
$$B_n:=A[X_1,\ldots,X_n]\to B,$$
and consequently a monomorphism
$$\Der_A(B)\to\Der_A(B_n)\cong B^n.$$
Denote by $\partial_j$ (the isomorphic images in $\Der_A(B_n)$ of) the 
generators of the free module $B^n$. For a given
$f\in{\cal B}$, choose $\delta$ and $\theta$ such that the condition to be in
${\cal B}$ holds, and express $\delta=\sum_jg_j\partial_j$ and 
$\theta=\sum_jh_j\partial_j$. Abusing slightly notation, we have
$$f(p)g_j(p)\,=\,h_j(p)$$
for all $p\in\pi^{-1}(x)$, by definition of $\delta$ and 
$\theta$. This means that $f$ is a rational function on $Y$, regular on fibers
$\pi^{-1}(x)$ for $x$ smooth. Therefore ${\cal B}\subset K(Y)$. ${\cal L}$ is a
${\cal B}$-module by setting $f\cdot\delta=\theta$, for the $\theta$ specified
through the definition of ${\cal B}$ by $f$ and $\delta$. Its annihilator is 
zero. On the other hand, ${\cal L}=\Der_A(B)$ is a finitely generated 
$B$-module if $A$ is noetherian. By a well-known criterion (cf \cite{Eise} 
p. $123$, cor. $4.6$), ${\cal B}$ is thus an integral extension of $B$. 

In conclusion, if $Y$ is normal, $B={\cal B}$ is characterized by Lie theoretic
methods, and the first part of the proposition follows.

If all fibers are smooth, ${\cal B}$ is even more readily identified with $B$.
If ${\cal L}$ is a free rank $1$ $B$-module with generator $\delta$ say, we 
have for $f\in{\cal B}$ by definition $f\delta=\theta$, but on the other hand
$f'\delta=\theta$ for some $f'\in B$. It follows $f=f'$ on all 
${\cal M}_{\phi}$, which permits to conclude in this case. \fin

\begin{rem}
In the smooth case, one can also argue in several different ways: in case
$Y=\Spec(B)$, $X=\Spec(A)$ and $\pi:Y\to X$ are smooth and $T\pi$ is 
surjective, the tangent bundles ${\cal T}_Y$ of $Y$, and ${\cal T}_X$ of $X$
fit into an exact sequence
$$0\to{\cal T}_{Y/X}\to {\cal T}_Y\to  \pi^*{\cal T}_X\to 0,$$
defining the relative tangent bundle ${\cal T}_{Y/X}$ whose regular sections 
are $\Der_A(B)$. One can show that $\Der_A(B)$ is an admissible $B$-Lie module
in the same way Amemiya showed it for Lie algebras defining pseudo foliations.

One can also show directly that $\Der_A(B)$ is ``strongly nowhere vanishing'' 
in the sense of Grabowski, in the smooth setting.
\end{rem}

Let us finish the {\it proof of theorem $1$:} it remains to show that two 
families of punctured curves $\pi:C\to S$ and $\pi':C'\to S$ with associated 
divisors $\triangle$ and $\triangle'$ such that 
$\phi:C\setminus\triangle\cong C'\setminus\triangle'$ as schemes over 
$\Spec(A)$ are isomorphic as schemes over $\Spec(A)$. 
Indeed, a point where $\phi$ is not defined is on some fiber of $\pi$ 
(the families are surjective), but there the restriction of $\phi$ defines 
a rational map of curves, thus there is no such point.  

\section{The stack of deformations}

It is well known that the moduli space of curves carries a structure of an 
algebraic {\it scheme} (the quasi-projective irreducible coarse moduli scheme),
or the structure of an {\it algebraic stack}. The first approach follows from
geometric invariant theory (GIT) and has the advantage that one stays within
the framework of standard algebraic geometry. The second approach stems from
Grothendieck's program of characterizing the functor of points 
$h_X:{\tt Sch}\to{\tt Sets}$, $S\mapsto\Hom_{\tt Sch}(S,X)$ associated to a 
fixed scheme $X$ within all functors $F:{\tt Sch}\to{\tt Sets}$. 
Some necessary 
conditions generalized from this example lead to the definition of a 
$\C$-stack, a device which generalizes schemes and which is important in the
study of classification problems in algebraic geometry where the 
``classifying space'' is not a scheme any more. Indeed, in the category of 
(algebraic) stacks, the moduli functor is representable, thus it does not only
posess a coarse, but a fine moduli stack. The disadvantage of the stack 
approach is that one has to (re)learn how to do (algebraic) geometry with 
stacks instead of schemes, varieties or manifolds. One advantage of the 
stack approach is that it includes (at least formally) the approach by 
coarse moduli spaces, cf rem. (3.19) \cite{LauMor}. The GIT is used in the
context of stacks in order to show that some stacks are algebraic.

Our understanding of the notion of algebraic or differentiable stacks is based
on \cite{Grot}, \cite{Metz}, \cite{Vist}.

When we address the question about some kind of structure of a variety 
on the space of deformations, it will lead naturally to the question whether 
the map $I$ respects these structures.

Let us define a stack of deformations of Lie algebras. It is some kind of
functor, more precisely a lax functor, which does not associate to an affine 
scheme $\Spec(A)$ the set of isomorphism classes of $A$-Lie algebras, 
but rather the groupoid of $A$-Lie algebras, keeping the data of the
iso- and automorphisms between/of $A$-Lie algebras. 
There are (at least) two points of view on stacks; let us use
here the one based on lax functors (or $2$-functors). 

As in any geometric theory, one needs first to fix a class of standard spaces
to which the spaces one wants to define, should be locally isomorphic. Led by
the example of the category of open sets of a manifold (the ones isomorphic to
$\R^n$ serving exactly the expressed need), one defines a Grothendieck 
(pre-)topology on a category ${\tt C}$ to be a collection of covering 
families $T(U)$ for each
object $U$ of ${\tt C}$ such that covering families contain isomorphisms, are
stable under base change, and are stable under refinement. The category of
interest for us will be $({\tt Aff}/\C)$, the category of affine schemes over 
$\C$, and as Grothendieck topologies on $({\tt Aff}/\C)$, we will regard the 
classical four topologies fpqc, fppf, \'etale and Zariski. The preceding 
abbreviations mean {\it fid\`element plat, quasi-compact}, and 
{\it fid\`element plat, de pr\'esentation finie} respectivement. This indicates
the conditions on a morphisms between affine schemes to be a member of the 
covering family. We refer to \cite{AltKlei} for basics about faithfully flat or
quasicompact morphisms, or morphisms of finite presentation. The four 
topologies are ordered from finest (fpqc) to coarsest (Zariski). See also
the precisions made by Kleiman and Vistoli on the fpqc topology in \cite{Vist}.

A {\it lax functor} associates to an affine scheme a groupoid. This 
association is not an ordinary functor, but a $2$-functor. Here the category
${\tt Aff}/\C$ is seen in a trivial way as a $2$-category (the  
$2$-morphisms are identities between compositions of $1$-morphisms), 
and the $2$-functor has values in the $2$-category of groupoids. 

More precisely, we define the {\it lax deformation functor} ${\cal D}{\rm ef}$
in the 
following way: let $U=\Spec(B)\in{\rm ob}\,({\tt Aff}/\C)$ be some affine 
scheme over $\C$. The lax functor ${\cal D}{\rm ef}:({\tt Aff}/\C)\to{\tt Gpd}$
associates to $U$ the groupoid ${\cal D}{\rm ef}(U)$ having as class of 
objects the $B$-Lie algebras and as morphisms only isomorphisms of $B$-Lie 
algebras. We forget for the moment that the Lie algebra should be isomorphic 
to a given one at a specified point. 
For a morphism $f:U'\to U$ in $({\tt Aff}/\C)$, we denote by 
$\tilde{f}:B\to B'$ the corresponding morphism of $\C$-algebras. To $f$, 
the lax functor ${\cal D}{\rm ef}$ associates a functor 
$f^*:{\cal D}{\rm ef}(U)\to{\cal D}{\rm ef}(U')$ given by 
$f^*({\mathfrak g})={\mathfrak g}\otimes_BB'$ (where the bracket on the tensor
product is the current algebra bracket) for any $B$-Lie algebra 
${\mathfrak g}$. Here $B'$ is seen as a $B$-algebra via $\tilde{f}$.

Let us verify the axioms of a lax functor ($2$-functor): 
let $g:U''\to U'$ be another 
morphism in $({\tt Aff}/\C)$ with associated morphism $\tilde{g}:B'\to B''$ of
$\C$-algebras. Then we have $g^*\circ f^*\simeq(f\circ g)^*$. Indeed, 
$(f\circ g)^*({\mathfrak g})={\mathfrak g}\otimes_BB''\simeq({\mathfrak g}
\otimes_BB')\otimes_{B'}B''=g^*(f^*({\mathfrak g}))$, where $B''$ (resp. $B'$)
is seen as a $B'$, $B$ (resp. $B$-) module via $\tilde{g}$, 
$\widetilde{(f\circ g)}=\tilde{g}\circ\tilde{f}$ (resp. $\tilde{f}$).

Now consider a third morphism $h:U'''\to U''$ ($\tilde{h}:B''\to B'''$). 
We have then a commutative diagram of $B'''$-Lie algebras

\vspace{.5cm}
\hspace{-.5cm}
\xymatrix{
h^*\circ g^*\circ f^*({\mathfrak g})=\left(({\mathfrak g}\otimes_BB')
\otimes_{B'}B''\right)\otimes_{B''}B'''  \ar[r] \ar[d] & ({\mathfrak g}
\otimes_BB'')\otimes_{B''}B'''  \ar[d] \\
(g\circ h)^*\circ f^*({\mathfrak g})=({\mathfrak g}\otimes_BB')
\otimes_{B'}B''' \ar[r]  & {\mathfrak g}\otimes_BB'''}
\vspace{.5cm} 

Let us show now that our lax functor is a $\C$-stack. It is here that the 
chosen Grothendieck topology on ${\tt Aff}/\C$ comes into play.
In order to be a $\C$-stack, the lax functor
should be on the one hand a sheaf of spaces, and on the other, one should be  
able to define in a unique way morphisms and objects from data given on a 
covering family (i.e. equivalence of descent categories). More precisely, 
the first condition is that for all $U\in{\rm ob}({\tt Aff}/\C)$, and all 
${\mathfrak g},{\mathfrak h}\in{\rm ob}({\cal D}{\rm ef}(U))$, the presheaf
$$\Hom({\mathfrak g},{\mathfrak h}):\,\,{\tt Aff}/\C\,\to{\tt Sets},$$
given by
$$(f:U'\to U)\,\mapsto\,\Hom_{{\cal D}{\rm ef}(U')}(f^*{\mathfrak g},
f^*{\mathfrak h}),$$
should be a sheaf, i.e. one should have an exact sequence of sets

\vspace{.5cm}
\xymatrix{
\Hom_{{\cal D}{\rm ef}(U)}({\mathfrak g},{\mathfrak h}) \ar[r] &
\Hom_{{\cal D}{\rm ef}(U')}(f^*{\mathfrak g},f^*{\mathfrak h}) \ar@<.5ex>[r]
\ar@<-.5ex>[r] &
\Hom_{{\cal D}{\rm ef}(U'')}(q^*{\mathfrak g},q^*{\mathfrak h}),}
\vspace{.5cm}
 
\noindent which lies over the ``exact'' sequence of objects in ${\tt Aff}/\C$

\vspace{.5cm}
\hspace{2cm}
\xymatrix{
U''=U'\times_{U}U' \ar@<.5ex>[r] \ar@<-.5ex>[r] &
U'  \ar[r] &
U,}
\vspace{.5cm}

\noindent or, equivalently, over the ``exact'' sequence of $\C$-algebras

\vspace{.5cm}
\hspace{2cm}
\xymatrix{
B   \ar[r]  &
B'  \ar@<.5ex>[r] \ar@<-.5ex>[r] &
B''=B'\otimes_BB',}
\vspace{.5cm}

$q$ being associated to the arrow $\tilde{q}:B\to B''$. Let us emphasize that
$f:U'\to U$ is here a covering family in the given Grothendieck topology, 
which means that the morphism $f$ has some geometric properties, but should be
moreover surjective (because one thinks $U'$ as $\coprod_i U_i'$, and the 
$U_i'$ should cover $U$).    

The second condition is that any descent data is {\it effective}, i.e. that 
for any covering family 
$f:U'\to U$, the category ${\cal D}{\rm ef}(U)$ is equivalent to the category 
of descent data in ${\cal D}{\rm ef}(U')$, i.e. to the category whose objects 
are the paires $(x',\phi)$ of an object $x'\in{\rm ob}({\cal D}{\rm ef}(U'))$ 
and an isomorphism 
$\phi:p_1^*x'\to p_2^*x'$, where $p_1,p_2:U''\to U'$ are the two projections,
such that $\phi$ satisfies the cocycle identity 
$p_{13}^*\phi=p_{23}^*\phi\circ p_{12}^*\phi$, where $p_{ij}$ are projections
$p_{ij}:U'''\to U''$ onto a choice of two factors in the triple fibered 
product. 

It is a non-trivial fact that the lax functor associating to an affine scheme
$U=\Spec(B)$ the category of $B$-modules is a $\C$-stack where ${\tt Aff}/\C$
carries the fpqc topology. This is the content of Grothendieck's theorem of 
{\it descente fid\`element plate} \cite{Grot}. One easily deduces (cf 
Expos\'e VIII, \S 2 in \cite{Grot}) the following

\begin{theo}
The deformation lax functor ${\cal D}{\rm ef}$ is a $\C$-stack.
\end{theo}

In this form, ${\cal D}{\rm ef}$ is certainly not an algebraic stack, nor a 
differentiable
stack in the sense of \cite{Metz}, as one needs a covering (algebraic) space.
Nevertheless, fixing an infinite dimensional Lie algebra $\Der_A(B)$ with
$X_0=\Spec(B)$, which is to be deformed, one can restrict ${\cal D}{\rm ef}$ 
to the $\C$-stack of deformations of $\Der_A(B)$ by considering first the full 
subcategory ${\tt Aff}_{X_0}/\C$ of affine schemes which have $X_0$ as fiber
over a closed point. Let $f:{\tt Aff}_{X_0}\to{\tt Aff}$ be the inclusion 
functor. Then one takes the direct image \cite{Gir} p. 83, \cite{Metz} p. 16
$f_*{\cal D}{\rm ef}$ as deformation lax functor,
and it is a standard matter to prove that $f_*{\cal D}{\rm ef}$ is still a 
$\C$-stack (in the induced Grothendieck topology on ${\tt Aff}_{X_0}$. 
This stack $f_*{\cal D}{\rm ef}$ is a more natural candidate for being an 
algebraic stack. It would be interesting to examine whether the miniversal 
deformation space \cite{FiaFuc} of a Lie algebra ${\mathfrak g}:=\Der_A(B)$ 
with finite dimensional $H^2({\mathfrak g},{\mathfrak g})$ could serve as 
such a a covering space in order to show that in this situation, 
$f_*{\cal D}{\rm ef}$ is an algebraic stack.  

\section{The morphism $I$}

In this section, we will show that the map $I$ gives a morphism of 
$\C$-stacks from the moduli stack ${\cal M}_{g,n}$ to the deformation stack 
${\cal D}{\rm ef}$, at least in any 
topology which is more coarse than the fppf topology (for example, in the 
\'etale topology, in which the moduli stack is usually considered).

\begin{theo}
The map $I$ defines a morphism of $\C$-stacks from ${\cal M}_{g,n}$ to 
${\cal D}{\rm ef}$ in any topology which is more coarse than the fppf topology.
\end{theo}

\pr Recall that the {\it moduli stack} ${\cal M}_{g,n}$, regarded as a lax 
functor, associates to an affine scheme $S$ over $\C$ the groupoid of proper 
flat families of curves $\pi:C\to S$ with sections 
$\sigma_1,\ldots,\sigma_n:S\to C$. The map $I$ is now regarded as morphism
of lax functors, i.e. for all affine schemes $U$,
$$I(U):{\cal M}_{g,n}(U)\,\to\,{\cal D}{\rm ef}(U)$$ 
is the functor associating to the family of curves $\pi:C\to U$ with sections 
$\sigma_1,\ldots,\sigma_n:U\to C$ the Lie algebra $\Der_A(B)$ such that 
$\Spec(B)=C\setminus\triangle$ with $\triangle=\bigcup_{i=1}^n\sigma_i(U)$ 
and $U=\Spec(A)$. 
Observe that $I(U)$ is a functor here, because ${\cal M}_{g,n}(U)$ and 
${\cal D}{\rm ef}(U)$ are groupoids. Now in order to be a morphisms of stacks, 
one needs furthermore a compatibility on $2$-morphisms, more precisely, for a
every covering family $f:U'\to U$ in ${\tt Aff}/\C$, 
there exists in the diagram

\vspace{.5cm}
\hspace{2.5cm}
\xymatrix{
{\cal M}_{g,n}(U)  \ar[r]^{I(U)} \ar[d]^{f^*} & {\cal D}{\rm ef}(U) 
\ar[d]^{f^*} \\
{\cal M}_{g,n}(U')  \ar[r]^{I(U')} & {\cal D}{\rm ef}(U')}
\vspace{.5cm}

a natural isomorphism $\alpha(f):f^*\circ I(U)\to I(U')\circ f^*$ between the 
two possible compositions of functors. This last statement is true in the 
fppf topology, because of lemma $1$, the compatibility of the fppf topology 
with base change, and the fact that 
$$(X\times_UU')\setminus\triangle'\cong(X\setminus\triangle)\times_UU'$$
which follows from the compatibility of the embedding $X\hookrightarrow \P^N$ 
with base change. This concludes the proof of theorem $3$.\fin

\begin{rem}
The morphism $I$ is almost a monomorphism; indeed, for integral normal $B$ 
we can apply Pursell-Shanks theory and use theorem $1$. But for a general 
affine base scheme $S=\Spec(A)$, one cannot hope to have monomorphy of $I$.
For example, a non-reduced scheme and its reduction may well have 
isomorphic tangent Lie algebra, cf remark on p. 313 in \cite{HauMue}.

Theorem $1$ shows however that $I$ is a monomorphism for the stack of families
of smooth curves.  
\end{rem}

\section{Lie algebra cohomology sheaves on moduli stack}

In this section, we will examine natural sheaves on the deformation 
stack ${\cal D}{\rm ef}$ arising from the cohomology of Lie algebras. The map
$I$ permits to pull them back to ${\cal M}_{g,n}$. 

Let ${\mathfrak g}$ be an $A$-Lie algebra, where $A$ is some $\C$-algebra.
Let $M$ be an $A$-module carrying an $A$-linear action 
of ${\mathfrak g}$. Then denote by $H_A^*({\mathfrak g};M)$ the cohomology of 
${\mathfrak g}$ with values in $M$, where cochains are supposed to be 
$A$-linear. There is a base change formula for this kind of cohomology space:

\begin{lem}
Let $A\to B$ be a morphism of commutative unital $\C$-algebras such that 
$B$ is an $A$-module of finite presentation. Then
$$H_B^*({\mathfrak g}\otimes_AB;M\otimes_AB)\,\cong\,H_A^*({\mathfrak g};M)
\otimes_AB.$$
\end{lem}

\pr Recall that the cochain complex computing $H_A^*({\mathfrak g};M)$ 
consists of spaces
$\Hom_A(\Lambda_A^p({\mathfrak g}),M)$, where $\Lambda_A^p({\mathfrak g})$ 
means the skewsymmetric tensor product over $A$ of $p$ factors of 
${\mathfrak g}$. The proof is easily deduced from three steps:
$$({\mathfrak g}\otimes_AB)\otimes_B({\mathfrak g}\otimes_AB)\,\cong\,
{\mathfrak g}\otimes_A{\mathfrak g}\otimes_AB,$$
$$\Hom_B({\mathfrak g}\otimes_A{\mathfrak g}\otimes_AB,X)
\,\cong\,
\Hom_A({\mathfrak g}\otimes_A{\mathfrak g},X),$$
where $X$ is a $B$-module, on the RHS viewed as an $A$-module, and
$$\Hom_A({\mathfrak g}\otimes_A{\mathfrak g},M\otimes_AB)\,\cong\,   
\Hom_A({\mathfrak g}\otimes_A{\mathfrak g},M)\otimes_AB,$$
where we used in the last step the finite presentation of $B$ as an $A$-module.
\fin

\begin{cor}
Let $B=\C$ and $A\to\C$ be a retraction of the unity $\C\to A$. Then
$$H_{\C}^*({\mathfrak g}\otimes_A\C;M\otimes_A\C)\,\cong\,
H_A^*({\mathfrak g};M)\otimes_A\C,$$
meaning that the fiber of the cohomology sheaf is the cohomology of the fiber
of the family of Lie algebras, cf lemma $2$.
\end{cor}

The above lemma means that the prescription 
$${\mathfrak g}\mapsto H_A^*({\mathfrak g};M)$$
for a fixed module $M$, defines a cartesian sheaf of vector spaces $H^*(M)$ 
on the stack ${\cal D}{\rm ef}$ in the fppf topology, and a cartesian sheaf 
$I^*(H^*(M))$ on the stack ${\cal M}_{g,n}$. The fixed module $M$ 
is supposed to be an 
$A$-module and a ${\mathfrak g}$-module for any $A$-Lie algebra 
${\mathfrak g}$. Examples are the ground field $M=\C$ (``trivial 
coefficients''), the Lie algebra itself $M={\mathfrak g}$ (``adjoint 
coefficients''), or algebraic constructions with these (like the symmetric
algebra $M=S^*{\mathfrak g}$ etc).

What is known about these cohomology spaces ?

It is known that $H^1(\C)=0$ on the locus of smooth families, because the Lie 
algebra of vector fields on a smooth curve is simple, see \cite{BFM} p. $21$.
I can also show that $H^2(Mer_k;\C)\cong H^1(\Sigma_k)$, 
the singular cohomology of the curve 
$\Sigma_k:=\Sigma\smallsetminus\{p_1,\ldots,p_k\}$
for any Lie algebra $Mer_k$ of meromorphic vector fields on the compact 
connected Riemann surface $\Sigma$ with possible p\^oles in the $k$ fixed 
points $p_1,\ldots,p_k$. The adjoint cohomology seems to be trivial.
All these are {\it algebraic} or {\it discrete} cohomology computations.
For {\it continuous cohomology}, i.e. cohomology with continuous cochains,
the Lie algebras of meromorphic vector fields $Mer_k$ are endowed with 
the subspace topology from the Lie algebra $Hol(\Sigma_k)$ of all holomorphic 
vector fields on $\Sigma_k$. The continuous cohomology of $Hol(\Sigma_k)$
(all degrees and values in various modules) is known by work of Kawazumi, 
but it seems artificial to endow these algebraic objects with topology.

One deduces from these results that the principle that cohomology should only
increase in isolated points is not true in this infinite dimensional setting.
Indeed, Fialowski and Schlichenmaier {\it loc. cit.} have a family of genus
$0$ curves with $3$ points which degenerates to the {\it Witt algebra} 
$$\Vect_{\rm pol}(S^1)\,=\,\bigoplus_{i\in\Z}\C x^{i+1}\frac{d}{dx}.$$
It is known that $H^2(\Vect_{\rm pol}(S^1);\C)$ is of dimension one, but on
the other hand, $H^2(Mer_3;\C)$ is of dimension two.  

\section{Deformations of Lie algebras arising from families of singular 
plane curves}

There is a way to incorporate families having as some of their members
 non stable curves such as cusps, for example, into the picture. 
Up to now, we took into account only
families of singular curves which constitute the boundary of the moduli stack.
The other approach is to normalize the family first in
a certain way and to map the thus obtained family of smooth curves into 
${\mathcal M}_{g,n}$. 
This corresponds
to the normalization which occurs in the procedure defined by Fialowski and
Schlichenmaier {\it loc. cit.}.

Consider families of irreducible plane curves with nodes and cusps as their 
only singularities. Arbitrary plane curves of degree $d$ are parametrized by a 
Hilbert scheme $\Hilb(d)$ which can be identified with $\P^N$ for $N=d(d+3)/2$.
Indeed, a curve is just represented by its homogeneous equation, up to scalars.
$\Hilb(d)$ is a fine moduli space and posesses a universal family, cf
\cite{HarMor}. 

In the following, we are inspired by the introduction of \cite{Gal}. 
Let $\Sigma^d_{k,n}\subset\P^N$ be the closure in the Zariski 
topology of the set of reduced and irreducible plane curves of degree $d$
with $k$ cusps and $n$ nodes. Let $\Sigma\subset\Sigma^d_{k,n}$ be an 
irreducible component of $\Sigma^d_{k,n}$. Denote by $\Sigma_0$ the open 
set of $\Sigma$ of points $s\in\Sigma$ such that $\Sigma$ is smooth at $s$
and such that $s$ corresponds to a reduced and irreducible plane curve of 
degree $d$ with $n$ nodes and $k$ cusps and no further singularities. 
The restriction of the universal family ${\mathcal S}_0\to\Sigma_0$ is an 
{\it equigeneric} family of curves, i.e. the genus $g$ of the normalizations
of the fibers of the family is constant. This property implies that one gets
a family of smooth curves by normalizing the total space, cf theorem (2.5) in
\cite{DiaHar}. In this way, there is a regular map of schemes from 
${\mathcal S}_0$ to the scheme of moduli of curves, and a morphism of algebraic
stacks ${\mathcal S}_0\to {\mathcal M}_g$, where ${\mathcal S}_0$ is seen as
an algebraic stack via its functor of points and $g=\left(\begin{array}{c}
d-1 \\ 2 \end{array}\right)-k-n$. One can introduce $m$ marked points in order 
to get a morphism of algebraic stacks
$$N:{\mathcal S}_0^m\to{\mathcal M}_{g,m}.$$
The composition 
$$I\circ N:{\mathcal S}_0^m\to{\mathcal M}_{g,m}\to{\mathcal D}ef$$ 
gives then the construction of Fialowski and 
Schlichenmaier in a conceptual way: the map $N$ consists of normalizing, 
i.e. desingularizing the singular fiber, and the map $I$ consists of extracting
the $m$ marked points from the family of smooth curves and then taking the Lie 
algebra of vector fields tangent to the fiber.
 
By its construction, it is clear that the same framework applies to flat 
equigeneric families of projective curves with all fibers reduced, where 
parameter and total space are reduced separated schemes of finite type 
over $\C$, see \cite{DiaHar} pp. 436--437.

\section{Examples}

This section is based on discussions with Martin Schlichenmaier.

\noindent(1)\quad Let us first consider a rather trivial example. 
Let ${\mathfrak g}$ be a non-abelian Lie algebra with Lie bracket $[,]$. 
Define a global deformation of ${\mathfrak g}$ over the affine line 
$A=\Spec(\C[t])$ by $({\mathfrak g}\otimes_{\C}\C[t],[,]_t)$ with the 
bracket $[,]_t=(t-1)[,]$.
This is a global deformation of ${\mathfrak g}$ which is non-trivial in the 
sense that there are two fibers which are non-isomorphic as $\C$-Lie algebras
(namely the fibers above $t=0$ and $t=1$; the Lie algebra above $t=1$ is 
abelian). But admittedly, it is trivial in some sense. 

\noindent(2)\quad Let $\pi:C\to S$ be a smooth family of curves with 
marked points $\sigma_1,\ldots,\sigma_n:S\to C$. We explained in sections $2$ 
and $3$ how to associate to $\pi:C\to S$ a global deformation ${\mathfrak g}_A$
of Lie algebras of the Lie algebra ${\mathfrak g}$ of regular vector fields on 
$C_0\setminus(\sigma_1(0)\cup\ldots\cup\sigma_n(0))$, with $\Spec(A)=S$. 
In case ${\mathfrak g}_A$ is a not a free $A$-module, this deformation cannot
be trivial (as a global deformation of Lie algebras).    

\noindent(3)\quad This is the example of Fialowski-Schlichenmaier 
{\it loc. cit.}. Elliptic curves can be parametrized in the complex 
projective plane $\P^2_{\C}$ by $e_1$, $e_2$ and $e_3$ such that the curve is 
given by the equation
$$Y^2Z\,=\,4(X-e_1Z)(X-e_2Z)(X-e_3Z)$$
with $e_1+e_2+e_3=0$ and $\triangle=16(e_1-e_2)^2(e_1-e_3)^2(e_2-e_3)^2\not=0$.
$\triangle\not=0$ assures that the curve is non-singular. These equations form 
a family of curves over $B:=\{(e_1,e_2,e_3)\,|\,e_1+e_2+e_3=0,\,\,\,\,
e_i\not=e_j\,\,\forall i\not=j\}$. 
Completing the parameter space $B$ to $\hat{B}:=
\{(e_1,e_2,e_3)\,|\,e_1+e_2+e_3=0\}$ admits singular cubics. Partial 
degeneration (e.g. $e_1=e_2\not=e_3$) leads to the {\it nodal cubic}
$$E_N\,:\,Y^2Z\,=\,4(X-eZ)^2(X+2eZ),$$
while overall degeneration (e.g. $e_1=e_2=e_3$) leads to the {\it cuspidal 
cubic}
$$E_C\,:\,Y^2Z\,=\,4X^3.$$
The nodal cubic is singular, but still stable and constitutes therefore a 
point of the boundary in the Deligne-Mumford compactification of the moduli 
space, while the cuspidal cubic is not. As the genus of the desingularization 
(i.e. normalization) decreases, the desingularization of $E_N$ and $E_C$ is
the projective line $\P^1_{\C}$. Fialowski and Schlichenmaier extract from the 
above family of curves a family over $\C[t]$ and compute explicitly the 
relations of the Lie algebra of vector fields on the family of elliptic curves.
It turns out that the generic fiber Lie algebra is not isomorphic to the Witt
algebra, while the fiber at $t=0$ is.

The point which makes this example less trivial than the preceeding is that
in any neighborhood of $0$, the restriction of the global deformation to this 
neighborhood remains non-trivial, while in the first two examples, there are 
neighborhoods where the deformation is trivial.

Fialowski and Schlichenmaier pointed out in {\it loc. cit.} that this is an
example of an infinitesimally rigid Lie algebra (because the Witt algebra has 
trivial second cohomology space with adjoint coefficients) which has 
non-trivial global deformations.     

The maps $I$ and $N$ permit in principle to define lots of examples of this 
third type. The hard part is to compute the relations of the so defined Lie
algebras in order to render the examples explicit.


\begin{thebibliography}{30}
\bibitem{AltKlei} A. Altman, S. Kleiman, Introduction to Grothendieck Duality 
Theory. Lecture Notes in Mathematics {\bf 146} Springer Verlag Heidelberg 1970
\bibitem{Amem} Ichiro Amemiya, Lie algebra of vector fields and complex 
structures. J. Math. Soc. Japan {\bf 27}, 4 (1975) 545--549
\bibitem{BFM} A. Beilinson, B. L. Feigin, B. Mazur, Algebraic conformal field 
theory. preprint
\bibitem{DiaHar} Steven Diaz, Joe Harris, Ideals of singular plane curves.
Trans. AMS {\bf 309}, 2 (1988) 433--468
\bibitem{Edid} Dan Edidin, Notes on the Construction of the Moduli Space of 
Curves.  Recent progress in intersection theory (Bologna, 1997),  85--113, 
Trends Math., Birkh\"auser Boston, Boston, MA, 2000
\bibitem{Eise} David Eisenbud, Commutative algebra. With a view towards 
algebraic geometry, Springer Graduate Texts in Mathematics 150, New York 1995
\bibitem{FiaFuc} Fialowski, Alice, Fuchs, Dimitry, Construction of Miniversal 
Deformations of Lie Algebras, J. Funct. Anal. {\bf 161} (1999) 76--110
\bibitem{FiaSch} Alice Fialowski, Martin Schlichenmaier, Global Deformations
of the Witt algebra of Krichever-Novikov type. Commun. Contemp. Math. Vol. 
{\bf 5}, No. 6 (2003) 921--945
\bibitem{Gal} Concettina Galati, Number of Moduli of irreducible families of 
plane curves with nodes and cusps. preprint on the arXiv: {\tt 0704.0618}
\bibitem{Grab} Janusz Grabowski, Isomorphisms and Ideals of the Lie algebra of
 Vector Fields. Invent. Math. {\bf 50} (1978) 13--33
\bibitem{Gir} Jean Giraud, Cohomologie non ab\'elienne. Die Grundlehren der
mathematischen Wissenschaften in Einzeldarstellungen, Band {\bf 179}, Springer
Verlag Heidelberg 1971 
\bibitem{Grot} Alexander Grothendieck et al., SGA 1: Rev\^etements \'etales
et Groupe Fondamental, Lecture Notes in Mathematics {\bf 224} Springer Verlag 
Heidelberg 1971
\bibitem{HarMor} Joe Harris, Ian Morrison, Moduli of Curves, Springer GTM 
{\bf 187}, 1998 
\bibitem{Hart} Robin Hartshorne, Algebraic Geometry. Springer Graduate Texts 
in Mathematics {\bf 52} (1977)
\bibitem{HauMue} Herwig Hauser, Gerd M\"uller, Affine varieties and Lie 
algebras of vector fields. Manuscripta Math. {\bf 80} (1993) 309--337
\bibitem{Knut} Finn F. Knudson, The projectivity of the moduli space of stable curves II: the stacks $M_{g,n}$. Math. Scand. {\bf 52} (1983) 161--199
\bibitem{LauMor} G\'erard Laumon, Laurent Moret-Bailly, Champs alg\'ebriques.
Ergebnisse der Mathematik und ihrer Grenzgebiete, Vol. {\bf 39}, Springer
Verlag Berlin 2000
\bibitem{Metz} David Metzler, Differentiable stacks. {\tt math.DG/0306176}
\bibitem{Sch} Martin Schlichenmaier, Oleg Sheinman, Knizhnik-Zalmolodchikov
equations for positive genus and Krichever-Novikov algebras. Russ.
Math. Surv. {\bf 59}, 4 (2004) 737--770
\bibitem{Sieb} Thomas Siebert, Lie algebras of derivations and affine algebraic
 geometry over fields of characteristic $0$. Math. Ann. {\bf 305} 
(1996) 271--286
\bibitem{Vist} Angelo Vistoli, Notes on Grothendieck topologies, fibered 
categories and descent theory. {\it Fundamental Algebraic Geometry: 
Grothendieck's FGA Explained}, ed. Barbara Fantechi et al., AMS 2005
\end{thebibliography}
\end{document}